\documentclass{amsproc}

\usepackage{amsmath,amscd,amssymb,amsfonts,latexsym,wasysym, mathrsfs, mathtools,hhline,color}
\usepackage[all, cmtip]{xy}

\usepackage{bm}

\newtheorem{theorem}{Theorem}[section]
\newtheorem{lemma}[theorem]{Lemma}

\theoremstyle{definition}
\newtheorem{definition}[theorem]{Definition}
\newtheorem{example}[theorem]{Example}
\newtheorem{que}[theorem]{Question}
\newtheorem{conj}[theorem]{Conjecture}
\newtheorem{prop}[theorem]{Proposition}
\newtheorem{cor}[theorem]{Corollary}
\newtheorem{assumption}[theorem]{Assumption}

\theoremstyle{remark}
\newtheorem{remark}[theorem]{Remark}

\numberwithin{equation}{section}

\def\be{\begin{equation}}
\def\ee{\end{equation}}

\def\bt{\begin{theorem}}
\def\et{\end{theorem}}

\def\bc{\begin{cor}}
\def\ec{\end{cor}}

\def\br{\begin{remark}}
\def\er{\end{remark}}

\def\bp{\begin{prop}}
\def\ep{\end{prop}}

\def\bl{\begin{lemma}}
\def\el{\end{lemma}}

\def\bex{\begin{example}}
\def\eex{\end{example}}

\def\bd{\begin{definition}}
\def\ed{\end{definition}}

\newcommand{\R}{\mathbb{R}}
\newcommand{\C}{\mathbb{C}}
\newcommand{\Z}{\mathbb{Z}}
\newcommand{\Hom}{\mathrm{Hom}}
\newcommand{\Q}{\mathbb{Q}}
\newcommand{\K}{\mathbb{K}}
\newcommand{\CP}{\mathbb{CP}}

\newcommand{\V}{\mathcal{V}}
\newcommand{\cha}{\mathrm{char}}

\newcommand{\p}{\overline{\mathbb{F}}_p}
\newcommand{\orb}{\mathrm{orb}}
\newcommand{\M}{\mathbf{M}}
\newcommand{\Tor}{\mathrm{Tor}}
\newcommand{\Pic}{\mathrm{Pic}}

\newcommand{\tor}{\mathrm{tor}}
\newcommand{\homo}{\mathrm{Hom}}
\newcommand{\sV}{\mathcal{V}}
\newcommand{\sL}{\mathcal{L}}

\newcommand{\sA}{\mathcal{A}}
\newcommand{\im}{\mathrm{im}}
\newcommand{\sO}{\mathcal{O}}

\newcommand{\rank}{\mathrm{rank }}
\newcommand{\coker}{\mathrm{coker}}

\newcommand{\T}{\mathbf{t}}

\begin{document}
\title[]{$L^2$-type invariants for complex smooth quasi-projective varieties -- a survey}

\author{Yongqiang Liu}
\address{Yongqiang Liu: Institute of Geometry and Physics, University of Science and Technology of China, 96 Jinzhai Road, Hefei Anhui 230026 China} 
\email{liuyq@ustc.edu.cn}

\thanks{The author  is  supported  by National Key Research and Development Project SQ2020YFA070080, NSFC grant No. 12001511, the Project of Stable Support for Youth Team in Basic Research Field CAS (YSBR-001), the starting grant from University of Science and Technology of China, the project ``Analysis and Geometry on Bundles" of Ministry of Science and Technology of the People's Republic of China and  Fundamental Research Funds for the Central Universities.
The author would like to thank the referee for  a careful reading of the manuscript and for several constructive suggestions. }


\date{May 20, 2024}


\begin{abstract}
Let $X$ be a complex smooth quasi-projective variety with an  epimorphism $\nu \colon \pi_1(X)\twoheadrightarrow \Z^n$.
 We survey recent developments about the asymptotic behaviour of Betti numbers with any field coefficients and the order of the torsion part of singular integral homology of finite abelian covers of $X$ associated to $\nu$, known as the $L^2$-type invariants. We give relations between $L^2$-type invariants, Alexander invariants and cohomology jump loci. 
When $\nu$ is orbifold effective, we give explicit formulas for $L^2$-invariants at homological degree one in terms of geometric information of $X$.
We also propose several related open questions for hyperplane arrangement complement.
\end{abstract}

\maketitle

\section{Introduction}

There is a general  principle  to consider a classical invariant of a finite CW complex and to define its analog for some  covering space. This leads to the $L^2$-type invariants. 
Atiyah  \cite{Ati76} introduced
the notion of $L^2$-Betti numbers in the context of a regular covering  of a closed
Riemannian manifold. 
After that,  there has been vast literature for the $L^2$-invariant theory, see \cite{Luc02}. The $L^2$-invariant theory has attracted considerable interest and has relations to many other fields, such as operator theory and Algebraic $K$-theory.   A particular important result is L\"uck's approximation theorem \cite{Luc94}, which states  that the $L^2$-Betti numbers of the universal cover of a finite polyhedron can be found as limits of normalised Betti numbers  of finite sheeted normal coverings.

\medskip

In this paper, we focus on a rather simple case of L\"uck's approximation problem, see e.g. L\"uck's survey papers \cite{Luc13,Luc16} on this topic. Let $X$ be a connected finite CW complex with a fixed epimorphism $\nu\colon\pi_1(X)\twoheadrightarrow \Z^n$. 
Consider  a descending chain of subgroups
$$\Z^n=\Gamma_0 \supseteq \Gamma_1 \supseteq \Gamma_2 \supseteq \cdots$$
such that $\Gamma_j$ is a sub-group in $\Z^n$ with finite index and $\bigcap_{j=0}^\infty\Gamma_j=\{0\}$.
Let $X_j$ denote the finite covering of $X$ associated to the composed epimorphism  $\nu $ with the quotient map $\Z^n \twoheadrightarrow \Z^n/\Gamma_j$. 
We study the asymptotic behaviour of 
the $i$-th Betti number with coefficients in a field coefficient $\K$
and the logarithm of the cardinality of the torsion subgroup of the $i$-th homology group with integral coefficients, i.e., the two sequences \begin{center}
    $ \big(\frac{\dim H_{i}(X_j,\mathbb{K})}{|\Z^n/\Gamma_j|}\big)_{j\geq 0} $
and
$\big(\frac{\log|H_{i}(X_j,\mathbb{Z})_{\mathrm{tor}}|}{|\Z^n/\Gamma_j|}\big)_{j\geq 0}.$
\end{center}   
Here $ |H_{i}(X_j,\mathbb{Z})_{\mathrm{tor}}|$ denotes the order of the torsion part of $H_{i}(X_j,\mathbb{Z}).$ 

\medskip
One can find the following approximation problem proposed by L\"uck in \cite{Luc13,Luc16}.
\begin{que} \label{que L2} (Approximation Problem) 

\begin{itemize}
    \item[(i)] Does the sequence converge?
    \item[(ii)] If yes, is the limit independent of the chain?
    \item[(iii)] If yes, what is the limit? 
\end{itemize}
\end{que}

Luckily for us Question \ref{que L2} has positive answers in our case as we explain now. 
 For a subgroup $\Gamma \subset \Z^n$ of finite index, we set
$$\langle\Gamma\rangle=\min \left\{ \sum_{j=1}^n x_j^2  \mid x=(x_1,\ldots,x_n)\in \Gamma, x\neq 0 \right\} $$
Let $X^{\nu,\Gamma}$ denote the covering space of $X$ associated to the compositions of $ \nu$ and the quotient map $\Z^n\to \Z^n/\Gamma$.
 Consider the following limits:
\[ \alpha_i(X^{\nu},\K)\coloneqq\lim_{\langle\Gamma\rangle\rightarrow\infty}\frac{\dim H_{i}(X^{\nu,\Gamma},\mathbb{K})}{|\Z^n/\Gamma|} \]
for any field coefficients $\K$ and
\[ \M_i(X^{\nu})\coloneqq\limsup_{\langle\Gamma\rangle\rightarrow\infty}\frac{\log|H_{i}(X^{\nu,\Gamma},\mathbb{Z})_{\mathrm{tor}}|}{|\Z^n/\Gamma|}. \]

It is known that the first limit always exists for any field coefficients and does not depend on the choice of chains, see \cite[Theorem 3.2]{Luc13}.
When $\cha(\K)=0$, the well known L\"{u}ck's  approximation theorem \cite{Luc94} shows that the first limit exists  and does not depend on the choice of chains. 
On the other hand, it was proved by Le \cite{Le}
that the upper limit $\M_*(X^\nu) $ exists and can be computed by the Mahler measure of the $i$-th integral Alexander polynomial.
When $n=1$, the upper limit in the definition of $\M_*(X^\nu) $ can be replaced by the ordinary limit.
\bigskip

The $L^2$-type invariants are particular useful for the following Singer-Hopf Conjecture.
\begin{conj}[Singer-Hopf Conjecture] \label{conj SH}
  Let $X$ be an aspherical closed manifold with real dimension $2N$. Then $$(-1)^N \chi(X)\geq 0.$$   
\end{conj}
A connected CW complex is said to be aspherical if its universal cover is
contractible. Conjecture \ref{conj SH} is still wide open even in the realm
of aspherical complex smooth projective varieties. This is particularly surprising as, in
the projective case, one can use a large variety of tools coming from algebraic
geometry approaches to such conjecture. For the recent developments in this direction, we refer to \cite{AW,AMW,CL,DL,LMW}. 

Inspired by these developments, we focus on the aforementioned $L^2$-type invariants in the context of complex smooth quasi-projective varieties. These $L^2$-type invariants have a deep connection with cohomology jump loci.
The cohomology jump loci of a complex smooth quasi-projective variety   have been  intensively studied in recent years. In particular, the beautiful structure theorem for cohomology jump loci put strong constraints for the homotopy type of complex smooth quasi-projective variety, see Section \ref{section structure theorem} for a brief introduction. We also refer the readers to Budur and Wang's survey paper \cite{BW17}  for a comprehensive background on this topic.

These results for cohomology jump loci can be used to prove many good properties about the $L^2$-type invariants. For example, as an application of structure theorem, we show that  the upper limit in the definition of $\M_*(X^\nu) $ can be replaced by the ordinary limit when $X$ is a smooth complex quasi-projective variety (see Corollary \ref{cor structure theorem}), which is not true for general topological space. This result is new to the author's knowledge.

Since the aforementioned limits for these  $L^2$-type invariants exist, we rephrase Question \ref{que L2}(iii) as follows.
\begin{que}
        How to compute these $L^2$-type invariants in terms of the algebraic geometric information of $X$, when $X$ is a smooth complex quasi-projective variety?
\end{que}
In the case of hyperplane arrangement complement, one can ask if these $L^2$-type invariants are combinatorially determined. 
 Based on L\"{u}ck's  approximation theorem \cite{Luc94}, a related question is whether the fundamental group of hyperplane arrangement complement is residually finite?

\bigskip

The purpose of this survey is to provide a gentle introduction to these concepts and
results, with an emphasis  on applications and open questions.
In section 2, we recall the definition of Alexander type invariants and cohomology jump loci and study their relations with $L^2$-type invariants.
In section 3, we focus on the case when $X$ is a complex smooth quasi-projective variety. In section 3.1, we recall the structure theorem about cohomology jump loci and its applications to $L^2$-type invariants. 
Section 3.2 is devoted to cohomology jump loci and $L^2$-type invariants at homological degree one. 
Section 4 is devoted to hyperplane arrangements. We propose many open questions regarding if these invariants are combinatorially determined. Section 5 is devoted for proofs.

\medskip

\textbf{Notations.}
\begin{enumerate}
\item $\K$ is assumed to be an algebraically closed filed in the rest of  the paper and $\cha(\K)$ denotes the characteristic of the field $\K$.
\item For any prime number $p\geq 2$, let $\p$ denote the algebraic closure of the finite field $\mathbb{F}_p $ with $p$ elements.
\end{enumerate}

\section{Alexander invariants and cohomology jump loci}
Let $X$ be a connected finite CW complex with an epimorphism $\nu\colon\pi_1(X)\twoheadrightarrow \Z^n$. 
We study $L^2$-type invariants associated to the pair $(X,\nu)$.

\subsection{Alexander type invariants}

Let $X^\nu$ denote the covering space of $X$ associated to the epimorphism $\nu\colon\pi_{1}(X)\twoheadrightarrow\mathbb{Z}^n$.
The group of covering transformations of the covering space  $X^{\nu}$ is isomorphic to $\Z^n$ and acts on it.   By choosing fixed lifts of the cells of $X$ to $X^{\nu}$, we obtain a free basis for the cellular chain complex  of $ X^\nu$ as   $\Z[\T^{\pm}]$-modules, where $\Z[\T^{\pm}]$ is the multi-variable Laurentiu polynomial ring $\Z[t^{\pm}_1,\ldots,t^{\pm}_n]=\Z[\Z^n]$.   So the cellular chain complex of $X^{\nu}$, $C_{*}(X^{\nu}, \Z)$, is a bounded complex of finitely generated free $\Z[\T^{\pm}]$-modules:
\be \label{chain compelx}
 \cdots  \to  C_{i+1}(X^\nu, \Z) \overset{\partial_{i}}{\to} C_i(X^\nu, \Z) \overset{\partial_{i-1}}{\to} C_{i-1}(X^\nu, \Z)  \overset{\partial_{i-2}}{\to}   \cdots \overset{\partial_0}{\to} C_0(X^\nu, \Z)  \to 0 .
\ee

\bd \label{def Alexander}
The {\it $i$-th integral  Alexander homology} $ H_{i}(X^{\nu},\Z)$ of the pair $(X,\nu)$ is by definition the 
$i$-th homology of  $C_*(X^\nu,\Z)$ as complex of $\Z[\T^{\pm}]$-modules.
Similarly, for any field $\K$, the {\it $i$-th Alexander homology with $\K$-coefficient} $ H_{i}(X^{\nu},\K)$ of  the pair $(X,\nu)$ is by definition the 
$i$-th homology of  $C_*(X^\nu,\Z)\otimes_\Z \K$ as complex of $\K[\T^{\pm}]$-modules, where $\K[\T^{\pm}]$ denotes the multi-variable Laurent polynomial ring $\K[t^{\pm}_1,\ldots,t^{\pm}_n].$
\ed 

With the above free basis for $C_*(X^\nu,\Z)$, $\partial_i$ can be written down as a matrix with entries in $\Z[\T^\pm]$. Note that $\Z[\T^\pm]$ is a Noetherian UFD.
Let $\Delta_i(X^\nu)$ denote the greatest common divisor of all non-zero $(\rank \partial_i\times \rank \partial_i)$-minors of $\partial_i$.
When $\partial_i=0$, $\Delta_i(X^\nu)=1$ by convention. 

\bd 
$ \Delta_i(X^\nu)$  is called {\it the $i$-th integral Alexander polynomial of the pair $(X,\nu)$}.  
 \ed
Note that $ \Delta_i(X^\nu)$ is defined uniquely up to  multiplication with a unit in $\Z[\T^\pm]$.

\br  The complex $C_*(X^\nu,\Z)$ determines $\alpha_*(X^\nu,\K)$ and $\M_*(X^\nu)$. 
In fact, for any $\Gamma \subset \Z^n$ with finite index, consider $\Z[\Z^n/\Gamma]$ as a quotient ring of $\Z[\T^\pm]$.
Then $H_*(X^{\nu,\Gamma}, \Z)$ can be computed by the following complex (viewed as a complex of abelian groups):
$$C_*(X^\nu, \Z) \otimes_{\Z[\T^\pm]}  \Z[\Z^n/\Gamma].$$
So it is no surprise that the $L^2$-type invariants have a deep connection with the Alexander invariants. 
\er

Next we recall the definition of Mahler measure.
Let $h\in \Z[\T^{\pm}]$ be a nonzero polynomial. The Mahler measure of $h$ is defined by
\[ \M(h)\coloneqq\int_{(S^1)^n} \log |h(s)|\mathrm{d}s,\]
where $\mathrm{d}s$ indicates the integration with respect to the normalized Haar measure, and $(S^1)^n$ is the multiplicative
subgroup of $n$-dimensional complex space $\C^n$ consisting of all vectors $(s_1, \ldots,s_n)$ with $|s_1|= \cdots =|s_n|=1$. 
When $n=1$ Mahler measure has the following equivalent definition. 
\bd \label{def MM} \cite{Mah60,Mah62}
For a non-zero polynomial $h(t)\in \Z[\T^\pm]$ with $n=1$, assume that \begin{center}
$h(t)=a_{d}t^{d}+a_{d-1}t^{d-1}+\cdots+a_{0}=a_{d}\cdot\prod_{j=1}^{d}(t-\alpha_{j})$ with $a_d\neq 0$, $a_0\neq 0$ and $\alpha_j\in \C$.
\end{center} 
We define the {\it Mahler measure} of $h$ to be
\[ \M(h)=\log(|a_{d}| \cdot \prod_{j=1}^{d}\max\{1,|\alpha_{j}|\}). \]
\ed

Recall that an element $h\in \Z[\T^\pm]$ is called a generalized cyclotomic polynomial \cite[page 47]{Sch} if it is of the form $h(t_1,\ldots,t_n)=  \Phi(\prod_{i=1}^n t_i^{a_i})$ (up to units in $\Z[\T^\pm]$), where $0\neq (a_1,\cdots,a_n)\in\Z^n$ and $\Phi$ is a cyclotomic polynomial in a single variable. 
The following theorem due to Schmidt characterizes when could the Mahler measure be zero. 
\bt\cite[Theorem 19.5]{Sch95}\label{thm cyclotomic}
For $h\in \Z[\T^\pm]$, $\M(h)=0$ if and only if $h$ is a product of generalized cyclotomic polynomials up to units.
\et

Based on the work of Schmidt \cite[Theorem 21.1]{Sch95}, Le  showed that \cite[Theorem 5]{Le} $\M_i(X^\nu)$ can be computed by the Mahler measure of $\Delta_i(X^\nu)$. 
When $n=1$, the following result is also proved by Raimbault \cite[Theorem 0.2]{Rai} independently.
\begin{theorem}\label{thm MM}  \cite[Theorem 5]{Le}
Let $X$ be a connected finite CW complex with an epimorphism $\nu\colon\pi_1(X)\twoheadrightarrow \Z^n$.
Then for any $i$
    we have $$\M_i(X^\nu)= \M(\Delta_i(X^\nu)).$$ 
   If $n=1$ then $\limsup$ in the definition of $\M_i(X^\nu)$ can be replaced by the ordinary $\lim$.
\end{theorem}

Inspired by Theorem \ref{thm cyclotomic}, we prove the following proposition. 
\bp \label{prop lim} Let $X$ be a connected finite CW complex with an epimorphism $\nu\colon\pi_1(X)\twoheadrightarrow \Z^n$. If $\Delta_i(X^\nu)$ is a product of some generalized cyclotomic polynomials with a non-zero integer $c_i$, then $\limsup$ in the definition of $\M_i(X^\nu)$ can be replaced by the ordinary $\lim$ and $\M_i(X^\nu)=\log \vert c_i \vert$.
\ep 
The proof of this proposition is given in section 5.

\subsection{Cohomology jump loci}
We first recall the definition of cohomology jump loci. For any algebraically closed field $\K$, the group of $\K$-valued  characters, $ \homo(\pi_1(X),\K^*)$, is a commutative affine algebraic group. Each character $\rho$ in $\homo(\pi_1(X),\K^*)$ defines a rank one local system on $X$, denoted by $L_{\rho}$. Since $\K^*$ is abelian, $\mathrm{Hom}(\pi_1(X),\K^*)$ in fact only depends on $H_1(X, \Z)$, the abelianization of $\pi_1(X)$ and $\homo(\pi_1(X),\K^*)=\homo(H_1(X,\Z),\K^*)$.
\bd
The cohomology jump loci of $X$ are defined as
$$\sV^i_k(X,\K)\coloneqq \lbrace \rho\in \homo(\pi_1(X),\K^*) \mid \dim_{\K} H^{i}(X, L_{\rho})\geq k \rbrace.$$
\ed

Cohomology jump loci  are closed sub-varieties of $\homo(\pi_1(X),\K^*)$ and homotopy invariants of $X$.
 In cohomological degree one, $\sV^1_k(X,\K)$  depends only on $\pi_1(X)$ (e.g. see \cite[Section 2.2]{Suc11}).

The map $\nu\colon \pi_1(X)\twoheadrightarrow  \Z^n$ induces an embedding $\nu^*\colon (\K^*)^n \hookrightarrow \homo(\pi_1(X),\K^*)$. For a tuple $\lambda=(\lambda_1,\ldots,\lambda_n)\in(\K^*)^n$, let $\nu^{-1} L_{\lambda}$ denote the corresponding rank one $\K$-local system on $X$ whose monodromy representation factors through $\nu$.

\bp  \label{prop betti number} 
With the  notations and assumptions  as above,  for any $i\geq 0$ and $\lambda \in (\K^*)^n$ being general we have 
$$ \alpha_i(X^\nu,\K)=\dim H_i(X, \nu^{-1} L_\lambda)=\rank H_i(X^\nu,\K).$$
Here $\rank H_i(X^\nu,\K)$ is the dimension  of the vector space $H_i(X^\nu,\K) \otimes \mathcal{F} $
over the fractional field $\mathcal{F}$ of $\K[\T^\pm]$. 
In particular, $\alpha_i(X^\nu,\K)$ is always an integer.  Moreover, if $\Delta_i(X^\nu)$ is not a unit in $\C[\T^\pm]$,  every irreducible component of the zero locus of $\Delta_i(X^\nu)$ in $(\C^*)^n$ is a  irreducible component of  $\sV^i_k(X,\C)  \cap (\C^*)^n$ (via $\nu^*$) with  codimension one for $k=\rank H_i(X^\nu,\C)+1.$
\ep
This proposition is not new and the claims are either already proved in \cite[Propsition 2.2]{LL23} or folklore.   We give a brief proof in section 5 due to the lack of precise references. For more relations between cohomology jump loci and Alexander type invariants, we refer the readers to Suciu's paper \cite{Suc11}.

\medskip

We also list some folklore properties of $L^2$-type invariants, which would be useful later.
\bp \label{prop folklore}   Let $X$ be a connected finite CW complex with a fixed epimorphism $\nu\colon\pi_1(X)\twoheadrightarrow \Z^n$. Then we have the following properties.
\begin{itemize}
    \item[(i)]   $\alpha_*(X^\nu,\K)$ and  $\M_*(X^\nu)$ are homotopy invariant of the pair $(X,\nu)$.
    \item[(ii)] {\it Euler number}: For any field $\K$, $$\sum_{i\geq 0} (-1)^i \alpha_i(X^\nu,\K) =\chi(X) .$$
    \item[(iii)] {\it Poincare duality}: If $X$ is  an oriented compact manifold with $\dim_\R X=N$, then   \begin{center}
            $\alpha_i(X^\nu,\K) =\alpha_{N-i}(X^\nu,\K)$ and $\M_i(X^\nu)=\M_{N-i-1}(X^\nu)$.
    \end{center}
    \item[(v)] {\it Zero-th $L^2$-type invariants}: $\alpha_0(X^\nu,\K)=0 $ and  $\M_0(X^\nu)=0$.
\end{itemize}
\ep 

\section{Complex smooth quasi-projective variety}
In this section, we always assume that $X$ is a complex smooth quasi-projective variety.
\subsection{Structure theorem}\label{section structure theorem}
The following  structure theorem for $\sV^i_k(X,\C)$ put strong constraints for the homotopy type of complex smooth quasi-projective variety. It is contributed by many people and we name a few here: Green and Lazarsfeld \cite{GL91}, Simpson \cite{Sim93}, Arapura \cite{Ara97}, Dimca, Papadima and Suciu \cite{DPS}, Artal Bartolo, Cogolludo-Agust\'in and Matei \cite{ACM13}, Dimca and Papadima \cite{DP14},  Schnell \cite{Sch}, etc. It is finalized by Budur and Wang in \cite{BW15,BW20}.
\bt \label{structure theorem} \cite{BW15,BW20} Let $X$ be a complex smooth quasi-projective variety. Then for any $i$ and $k$, $ \sV^i_k(X,\C)$ is a finite union of torsion translated sub-tori of $\Hom(\pi_1(X),\C^*)$. Moreover, for every positive dimensional irreducible component $V$ of $ \sV^i_k(X,\C)$, there exists a torsion point $\rho \in \homo(\pi_1(X),\C^*)$ and an algebraic map $f\colon X \to G$ with $G$ being a complex semi-abelian variety such that
$$ V= \rho \cdot \im\{ f^*\colon \homo(\pi_1(G),\C^*)\hookrightarrow \homo(\pi_1(X),\C^*)\}.$$
If $X$ is a smooth projective variety,  $G$ has to be an abelian variety, hence every irreducible component of $ \sV^i_k(X,\C)$ is even dimensional. 
\et
As a direct application of the above structure theorem, Theorem \ref{thm MM}, Proposition \ref{prop lim} and Proposition \ref{prop betti number}, we have the following results.
\bc  \label{cor structure theorem}
Let $X$ be a complex smooth quasi-projective variety with an epimorphism $\nu\colon\pi_1(X) \to \Z^n$.
Then for any $i$, $\Delta_i(X^\nu)$ is product of some generalized cyclotomic polynomials and a non-zero integer $c_i$. Then $\limsup$ in the definition of $\M_i(X^\nu)$ can be replaced by the ordinary $\lim$ and 
we have \begin{center}
    $\M_i(X^\nu)=\log \vert c_i \vert $. 
\end{center}
\ec 

Let $X$ be a complex smooth projective variety.
For any subgroup $\Gamma \subset \Z^n$ with finite index, $X^{\nu,\Gamma}$ is also a complex smooth projective variety. In particular, $H^i(X^{\nu,\Gamma},\Q)$ has pure Hodge structures. Instead of studying  the  asymptotic behaviour of
Betti number of $H^i(X^{\nu,\Gamma},\Q)$, one can also study the asymptotic behaviour of the Hodge numbers of $H^i(X^{\nu,\Gamma},\Q)$. 
 Let $h^{p,q}(X^{\nu,\Gamma})$ denotes the dimension of $ H^q(X^{\nu,\Gamma},\Omega^p_{X^{\nu,\Gamma}}) $. 
\bp \label{prop lim hodge}
With the above notations and assumptions,  for any $(p,q)$ we have
$$\lim_{\langle\Gamma\rangle\rightarrow\infty}\frac{h^{p,q}(X^{\nu,\Gamma})}{|\Z^n/\Gamma|}=\dim H^q(X,\Omega^p_X\otimes \sL'),$$
where $\sL'=L'\otimes \sO_X$ and $L'$ is a generic unitary rank one local system in $\im \nu^*$. Here $\nu^*$ is the embedding $(\C^*)^n \hookrightarrow \homo(\pi_1(X),\C^*)$ induced by $\nu$. 
\ep 

The proof of this proposition is given in section 5.
When  $\nu$  is the natural map to $H_1(X, \Z)$ modulo
torsions and the Albanese map of $X$ is semi-small,
 Di Cerbo and Lombardi \cite{CL} give an effective estimate for the asymptotic behaviour of the Hodge numbers and showed that $$\lim_{\langle\Gamma\rangle\rightarrow\infty}\frac{h^{p,q}(X^{\nu,\Gamma})}{|\Z^n/\Gamma|}=\begin{cases}
(-1)^q\chi(X,\Omega^p_X), & \mathrm{if}\ p+q=\dim_\C X,\\
0 , & \mathrm{else}.    
\end{cases} $$
Using the generic vanishing results proved by Popa and Schnell \cite{PS}, Di Cerbo and Lombardi's results are compatible with 
Proposition \ref{prop lim hodge}.

One may ask how about Proposition \ref{prop lim hodge} for  smooth quasi-projective varieties.  In this case,   $H^i(X^{\nu,\Gamma},\Q)$ has mixed Hodge structures.
 \begin{que} \label{que Hodge}
     Doe the limit for the Hodge numbers of  $H^i(X^{\nu,\Gamma},\Q)$ exist?
 \end{que}

Budur proved \cite[Theorem 1.8, Theorem 1.9]{Bu} the polynomial
periodicity for certain Hodge numbers of congruence covers of complex smooth (quasi-)projective variety based on  a structure theorem for line bundles \cite[Theorem 1.3, Theorem 1.4]{Bu}. 
His work shed light on the above question at least for the congruence covers.
For more references about polynomial
periodicity of congruence covers, we refer the readers to Suciu's survey paper \cite[Section 5]{Suc01}.

\subsection{Degree one case}
\subsubsection{Cohomology jump loci}
If we focus on cohomological degree one  jump loci, a  celebrity result  due to Beauville \cite{Bea}, Arapura \cite{Ara97} and Artal Bartolo, Cogolludo-Agust\'in and Matei \cite{ACM13} puts  even stronger constraints for its fundamental group.
To explain their results, we introduce the notion of orbifold maps. 
\bd Let $X$ be a smooth complex quasi-projective variety. An algebraic map $f\colon X \to C_{g,r}$  is called an orbifold map, if $f$ is surjective, has connected generic fiber and $C_{g,r}$ is a smooth algebraic curve of genus $g$ with $r$ points removed. We always assume that $C_{g,r} \neq \mathbb{CP}^1, \C$. There exists a maximal Zariski open subset $U\subset C_{g,r}$ such that $f$ is a fibration over $U$. Say $B=C_{g,r}-U$ (could be empty) has $s $ points, denoted by $\{q_1,\ldots,q_s\}$.  We assign
the multiplicity $\mu_j$ of the fiber   $f^{*}(q_j)$  (the $\gcd$ of the coefficients of the divisor $f^* q_j$) to the point $q_j$.   
Such orbifold map $f$ is called of type $(g,r,{\bm \mu})$, 
where ${\bm \mu}=(\mu_1,\ldots,\mu_s)$. When $B=\emptyset$, $\prod_{j=1}^s \mu_j=1$ by convention. 
 \ed
When $C_{g,r}$ is clear in the context, we simply write $C$.
Our definition of orbifold maps is a bit of different from the classical one, e.g. see \cite[Section 4.1]{DS14}. The classical one will be called the hyperbolic orbifold map, where we adopt the name from \cite{Delz}.
\bd An orbifold map $f\colon X \to C$ of   type $(g,r,{\bm \mu})$ is called hyperbolic, if 
$$2g+r-2+ \sum_{j=1}^s (1-\dfrac{1}{\mu_j})>0. $$
\ed 
This definition is used to rule out the case where $C$ is either $\C^*$ or Elliptic curve and $f$ has no multiple fibers.

The orbifold group $\pi_1^{\orb}(C_{g,r},{\bm \mu})$ associated to these data is  defined as
\[ \pi_{1}^{\orb}(C_{g,r},{\bm \mu})\coloneqq \pi_{1}(C_{g,r}\backslash \{q_{1}, \ldots, q_{s}\})/ \langle \gamma_{j}^{\mu_{j}}=1 \text{ for all } 1\leq j \leq s\rangle, \]
where $\gamma_{j}$ is a meridian of $q_{j}$.
An orbifold map  $f\colon X \to C$ of type $(g,r,{\bm \mu})$ induces an surjective map to the orbifold group (see \cite[Proposition 1.4]{ACM10})
$$ f_* \colon \pi_1(X) \twoheadrightarrow \pi_1^{\orb}(C_{g,r},{\bm \mu}).$$
Hence it induces an embedding for any $k\geq 1$
$$\sV^1_k(\pi_{1}^{orb}(C_{g,r},{\bm \mu}),\K) \to \sV^1_k(X,\K) .$$

\bigskip

The following theorem shows that every positive-dimensional component of $\sV^1_1(X,\C)$ is characterized by the orbifold map in this way.
The idea first appeared in Beauville’s work \cite{Bea} for the projective case and was extended to the
quasi-projective case by Arapura \cite{Ara97} and Artal Bartolo, Cogolludo-Agust\'in and Matei \cite{ACM13}. Further properties were found by Dimca \cite{Dim07,Dim09}, Dimca-Papadima-Suciu \cite{DPS}, etc.    
\bt \label{thm Ara} \cite{Ara97,ACM13,Dim07} 
Let $X$ be a complex smooth quasi-projective variety.   
\begin{itemize}
\item[(a)] We have
\[
\V^1_1(X,\C)=\bigcup_{f }{f^*} \sV^1_1(\pi_1^{\orb}(C_{g,r},{\bm \mu})),\C)\cup Z,
\]
where $Z$ is a finite set of torsion points and the union runs over all hyperbolic orbifold maps $f$.
\item[(b)] Fix an orbifold map $f\colon X\to C$ of type $(g,r,{\bm \mu})$. Then  for any  $\rho \in \homo(\pi_1^{\orb}(C_{g,r},{\bm \mu}),\C^*)$ we have
$$\dim H^1(X,f^* L_\rho)\geq \dim H^1(\pi_1^{\orb}(C_{g,r},{\bm \mu}),L_\rho)$$
and equality holds with finitely many (torsion) exceptions. 
\end{itemize}
\et
\br \label{rem complete C}
This theorem indeed gives a complete description of all positive-dimensional components contained in $\sV^1_k(X,\C)$ for all $k\geq 1$.
\er

When $\cha(\K)=p>0$, not much is known for $\sV^i_k(X,\K)$. So far all the known results concerns about cohomological degree one  jump loci. 
Delzant \cite{Delz} enabled  to generalize Theorem \ref{thm Ara}(a) to any algebraically closed $\K$-coefficient case when $X$ is a complex smooth projective variety.
\bt \label{thm Delz} \cite{Delz} Let $X$ be a complex smooth projective variety and $\K$ be an algebraically closed field. Then   we have
\[
\V^1_1(X,\K)=\bigcup_{f }{f^*} \sV^1_1(\pi_1^{\orb}(\Sigma_{g,r},{\bm \mu})),\K)\cup Z',
\]
where $Z'$ is a finite set of torsion points (depends on $\K$) and the union runs over all hyperbolic orbifold maps $f$. In particular, the structure theorem holds for $\sV^1_1(X,\K)$.
 \et 

The following theorem due to Li and the author  generalizes  Theorem \ref{thm Ara}(b) to any algebraically closed field coefficients. 

\bt \label{thm LL jump loci} \cite[Theorem 1.6]{LL21} Let $X$ be a complex smooth quasi-projective variety and $\K$ be an algebraically closed field. Given an orbifold map $f\colon X\to C$ of type $(g,r,{\bm \mu})$ and any  $\rho \in \homo(\pi_1^{\orb}(C_{g,r},{\bm \mu}),\K^*)$,
we have 
$$  \dim H^1(X,f^* L_\rho)\geq \dim H^1(\pi_1^{\orb}(C_{g,r},{\bm \mu}),L_\rho)$$
and equality holds with finitely many (torsion) exceptions. 
\et 

\br As like in Remark \ref{rem complete C}, Theorem \ref{thm Delz} and Theorem \ref{thm LL jump loci} together give a complete 
description for all positive-dimensional components contained in $\sV^1_k(X,\K)$ for all $k$ when $X$ is smooth projective. Hence the structure theorem holds for  $\sV^1_k(X,\K) $ for all $k\geq 1$.
\er


\subsubsection{BNS-invariants}
One may wonder how  about Theorem \ref{thm Delz} for  smooth quasi-projective varieties.
Following Delzant's approach in \cite{Delz,Delz10}, one needs  to study the BNS-invariants of $X$. Let us recall its definition.

Pick a finite generating set of $\pi_1(X)$ and let $\Omega(X)$ be the
corresponding Cayley graph of $\pi_1(X)$. Given an additive real character $\chi \colon \pi_1(X) \to \R$, let $\Omega_\chi(X)$ be the full subgraph on vertex set $\{g\in \pi_1(X)\mid  \chi(g)\geq 0\}$. Bieri, Neumann, and Strebel defined \cite{BNS}
a subset $\Sigma^1(X)$ of $H^1(X, \R) \setminus \{0\}$, nowadays called the BNS-invariant of $X$ as follows, see e.g. \cite[Definition 7.1]{PS10}.
\bd 
 The set $\Sigma^1(X)$ consists of those non-zero homomorphisms $\chi \colon \pi_1(X) \to \R$ 
for which the graph $\Omega_\chi(X)$ is connected.
\ed 

Papadima and Suciu studied the relations between $\Sigma^1(X)$ and the orbifold maps in \cite[Section 16]{PS10}, and the results are refined by Suciu recently via tropical method \cite{Suc23}. 
\bt \label{thm PS} Let $X$ be a complex smooth quasi-projective variety.
 Then $$\Sigma^1(X) \subset \big( \bigcup_f \im(f^*\colon H^1(C,\R)\to H^1(X,\R) ) \big)^c,$$
where the finite union runs over all hyperbolic orbifold maps $f$.
\et
In fact, a result of Delzant \cite{Delz10} shows that the above inclusion holds as an equality when $X$ is a smooth projective variety (or even a compact K\"ahler manifold). 
Based on a deep theorem due to Simpson \cite{Sim93B},  Delzant gave a complete description of BNS-invariants of complex compact K\"ahler manifold \cite{Delz10}, which are characterized by the hyperbolic orbifold maps.  Using this  description of BNS-invariants by orbifold maps, Delzant managed to prove Theorem \ref{thm Delz} in \cite{Delz}.
 \begin{que} \label{que 3} Inspired by these work, one may ask the following questions:
\begin{itemize}
  \item[(i)] Does \cite[Theorem 1]{Sim93B} hold for complex smooth quasi-projective varieties?
    \item[(ii)] Does the inclusion in Theorem \ref{thm PS}  hold as an equality?
    \item[(iii)] Does Theorem \ref{thm Delz} hold for complex smooth quasi-projective varieties?
\end{itemize}
\end{que}
Ideally a positive answer to question (i) (resp. (ii)) would imply a positive answer to question (ii) (resp. (iii)). 
The question (ii) was first raised by Papadima and Suciu in \cite{PS10}  and refined by  Suciu in \cite[Question 13.7]{Suc23} for hyperplane arrangement complements.

\subsubsection{$L^2$-type invariants}
To compute the $L^2$-type invariants of the pair $(X,\nu)$ at homological degree one, we need to connect $\nu$ with the orbifold maps.
\bd \label{def orbifold effective}
Let $X$ be a smooth complex quasi-projective variety with an epimorphism $\nu\colon \pi_1(X)\twoheadrightarrow \Z^n$. We say that $\nu$ is orbifold effective if there is an orbifold map $f\colon X\rightarrow C$ of type $(g,r,{\bm \mu})$ such that $\nu$ factors through $f_*$ as follows:
\[\xymatrix{
\pi_1(X)\ar@{->>}"1,3"^{\nu}\ar@{->>}[dr]_{f_*} & & \Z^n \\
 & \pi_1^{\orb}(\Sigma_{g,r}, {\bm \mu})\ar@{->>}[ur] &
}.\]
 We  say that $\nu$ is  orbifold effective by $f$ and call $\nu$ being of type $(g,r,{\bm \mu})$.
\ed
\br 
  As shown in \cite[Remark 4.3]{LL21}, if $\nu$ is orbifold effective by two different orbifold maps $f$ and $g$, then $f$ and $g$ must have the same type  $(g,r,{\bm \mu})$. So it is well defined to call $\nu$ being of type $(g,r,{\bm \mu})$.
\er
The following theorem by Li and the author  gave concrete formulas for the $L^2$-type invariants of $X$ at homological degree one when $\nu$ is orbifold effective.

\bt\label{main thm} \cite[Theorem 1.10]{LL21}, \cite[Theorem 1.3]{LL23}
Let $X$ be a complex smooth quasi-projective variety with an epimorphism $\nu\colon\pi_1(X)\twoheadrightarrow \Z^n$. Suppose that $\nu$ is  orbifold effective of type $(g,r, {\bm \mu})$. Let $\K$ be a field with $\cha(\K)=p\geq 0$. 
Then
$$  \alpha_1(X^\nu,\K)=  2g+r-2+ \#\{ j  \mid p\ \mathrm{ divides }\ \mu_j\} $$
and
$$  \M_1(X^\nu)=\sum_{j=1}^s \log \mu_j. $$
On the other hand, assume that either $X$ is projective or $H^1(X,\Q)$ has pure Hodge structure of type $(1,1)$. Then $\alpha_1(X^\nu,\K)>0$ or $  \M_1(X^\nu)>0$ happens only if $\nu$ is orbifold effective. 
\et

\br A positive answer to Question \ref{que 3}(iii) would enable us to drop the assumption that either $X$ is  projective or $H^1(X,\Q)$ has pure Hodge structure of type $(1,1)$ in the last claim of the above theorem.
\er 

\br Due to  Proposition \ref{prop betti number} and Corollary \ref{cor structure theorem}, the above theorem also implies the corresponding formulas for 
the rank of the Alexander homology and 
the leading coefficient of the integral Alexander polynomial at homological degree one.
\er 
  The above theorem is just about  $L^2$-type invariants at homological degree one. How about  higher homological degrees?
Following the structure theorem, we propose the following higher homological degree analogues. Let $V$ be a positive-dimensional irreducible component in $\sV^i(X,\C)$. According to the structure theorem there exists an algebraic map $f\colon X \to G$ with $G$ being a complex semi-abelian variety such that $V$ is contained in a translate of
$$  \im\{ f^*\colon\homo(\pi_1(G),\C^*)\hookrightarrow \homo(\pi_1(X),\C^*)\}.$$
\begin{que} \label{que higher degree} Note that $\pi_1(G)$ is  a finitely generated free abelian group. Consider the epimorphism $\nu=f_*\colon \pi_1(X)\to \pi_1(G)$ induced by $f$. Can one compute $\alpha_i(X^\nu,\K)$ and   $\M_i(X^\nu)$ by the geometric information of $f$? One may need to take the Stein factorization of $f$.
\end{que}

Using Theorem \ref{main thm} and Proposition \ref{prop folklore}, we give a complete description of the $L^2$-type invariants when $X$ is a complex smooth projective surface.
\bp Let $X$ be a complex smooth projective surface with an epimorphism $\nu\colon\pi_1(X)\twoheadrightarrow \Z^n$. 
Let $\K$ be a field with $\cha(\K)=p\geq 0$. If $\nu$ is not orbifold effective,
then  $$\alpha_i(X^\nu,\K) =\begin{cases}
\chi(X), & \mathrm{if}\ i=2,\\
 0, & \mathrm{else}
\end{cases} $$   and \begin{center}
  $\M_i(X^\nu)=0$ for all $i$.  
\end{center}
On other hand,  if $\nu$ is   orbifold effective, say of type $(g,r, {\bm \mu})$, then
 $$\alpha_i(X^\nu,\K) =\begin{cases}
2g+r-2+ \#\{ j  \mid p\ \mathrm{ divides }\ \mu_j\}, & \mathrm{if}\ i=1,3,\\
\chi(X)+4g+2r-4+ 2\#\{ j  \mid p\ \mathrm{ divides }\ \mu_j\}, & \mathrm{if}\  i=2,\\
 0 & \mathrm{ else}
\end{cases} $$   and
\begin{center}
    $\M_i(X^\nu) =\begin{cases}
\sum_{j=1}^s \log \mu_j, & \mathrm{if}\ i=1,2,\\
0 , & \mathrm{else}.
\end{cases} $
\end{center}
\ep

\section{Hyperplane arrangement} In this section, we always assume that  $X$ is the complement of a hyperplane arrangement $\sA$ in $\mathbb{CP}^N$ for some $N$ and we focus on a fixed epimorphism $\nu\colon\pi_{1}(X)\rightarrow\Z$ for $n=1$.

For any field $\K$, we denote $\nu_{\Z}$ and $\nu_\K$ the corresponding element in $H^{1}(X,\Z)$ and $H^1(X,\K)$, respectively. Note that $\nu_{\Z}\cup v_{\Z}=0,$ hence   $\nu_{\K}\cup \nu_{\K}=0$. 
Then we get two chain complexes by cup product
\[ \xymatrix{
(H^{\ast}(X,  \Z) , \cdot\nu_{\Z})\colon & H^{0}(X,\Z)\ar[r]^{\nu_{\Z}} & H^{1}(X,\Z)\ar[r]^{\nu_{\Z}} & H^{2}(X,\Z)\ar[r] & \cdots \\
(H^{\ast}(X,\K),\cdot \nu_{\K})\colon & H^{0}(X,\K)\ar[r]^{\nu_{\K}} & H^{1}(X,\K)\ar[r]^{\nu_{\K}} & H^{2}(X,\K)\ar[r] & \cdots
} \]
which are  called {\it the Aomoto complex} of $H^{\ast}(X,\Z)$ and $H^{\ast}(X,\K)$ associated to $\nu_{\Z}$ and { $\nu_{\K}$, respectively.

\bd 
With the above notations, we define {\it  the $i$-th Aomoto Betti number with $\K$-coefficients} as
\[ \beta_{i}(X,\nu_{\K})\coloneqq\dim_{\K}H^{i}(H^{\ast}(X,\K), \cdot\nu_{\K}). \]
and {\it  the order of the $i$-th Aomoto torsion} as
\[ \tau_i(X,\nu_\Z)\coloneqq | H^{i+1}(H^*(X,\Z), \cdot\nu_\Z)_{\tor}|. \]
Here the shift by 1 is due to the Universal Coefficient Theorem.
\ed

As shown by Orlik and Solomon \cite{OS} the cohomology ring $H^*(X,\Z)$ are determined by the combinatorics data of the hyperplane arrangement. 
Then $\beta_*(X,\nu_\K)$ and $\tau_*(X,\nu_\Z)$ are combinatorially determined,    once $\nu$, considered as an element in $H^1(X,\Z)$, is fixed. 
Papadima and Suciu discovered the following results in \cite{PS10}.
\bt \label{thm inequality} With the above notations and assumptions, for any field $\K$ we have $$  \alpha_i(X^\nu,\K) \leq \beta_i(X,\nu_\K).$$ 
Moreover, when $\K=\C$, the above inequality holds as a equality.
\et 
Papadima and Suciu \cite{PS10} indeed proved the above inequality for any finite CW complex. When $\K=\C$, the equality
was first proved (by different methods) by Cohen and Suciu \cite{CS}, Libgober \cite{Lib}, and
Libgober and Yuzvinsky \cite{LY} for $i=1$ and was generalized to the higher-degrees 
by Cohen and Orlik  \cite[Theorem 3.7]{CO}.

On the other hand, the author jointly with Li gave the following combinatorial upper bound for  $\M_*(X^\nu)$.
\bp \label{prop comb} \cite[Proposition 1.12]{LL21} Let $X$ be the complement of a hyperplane arrangement.
For any epimorphism $\nu\colon \pi_1(X) \to \Z$ and any $i\geq 0$, we have 
 $$ \exp(\M_i(X^\nu)) \mid \tau_i(X, \nu_\Z).$$ In particular, if $\nu$ is orbifold effective of type $(g,r,{\bm \mu})$, then \be \label{upper bound for mu}
  \prod_{j=1}^s \mu_j \mid \tau_1(X,\nu_\Z). \ee
\ep 
The proof of this proposition is based on Papadima and Suciu's work \cite[Theorem 12.6]{PS10} and two facts about hyperplane arrangement complement:
\begin{itemize}
    \item the hyperplane arrangement complement $X$ is  homotopy equivalent to a minimal CW complex, see \cite{DP03,Ran}.
    \item  $\alpha_i(X^\nu,\C)=\beta_i(X,\nu_\C)$ for any $i$ as in Theorem \ref{thm inequality}.
\end{itemize}

There are examples of hyperplane arrangements where $\alpha_i(X,\K)\neq \beta_i(X,\K)$ and $\prod_{j=1}^s\mu_j\neq\tau_1(X,\nu_\Z)$, see  \cite[Example 5.3]{LL21} and the next example.

\bex \label{ex 2}
Let $\mathcal{B}$ be the deleted $B_3$-arrangement in $\CP^2$ with defining equations $$zxy(x-y)(x-z)(y-z)(x-y-z)(x-y+z).$$ Order the
hyperplanes as the factors of the defining polynomial. Let $X$ be the complement of the arrangement.  
It was first discovered by Suciu \cite{Suc02} that $\V^1_1(X,\C)$  has a translated component 
$$V\coloneqq \rho \otimes \{t, t^{-1},1,t^{-1},t,t^2,t^{-2} \mid t\in \C^*\} $$
with $\rho = (1,-1,-1,-1,1,1,1)\in (\C^*)^7$, e.g. see \cite[Example 6.16]{Dim17}. Here we take $\{z=0\}$ as the hyperplane at infinity.
This component is induced by the orbifold map $f\colon X\rightarrow \C^*$ as following
\[ f([x,y,z])=\frac{x(y-z)(x-y-z)^2}{y(x-z)(x-y+z)^2} \]
and $f$ is  of type $(0,2,2)$. 
 Consider $\nu=(1,-1,0,-1,1,2,-2) \in H^1(X,\Z)$ induced by $f$. Then
\begin{center}
  $\alpha_1(X^\nu,\K) =\begin{cases}
1, & \mathrm{if}\ \cha(\K)=2,\\
0 , & \mathrm{if}\ \cha(\K)\neq 2,
\end{cases} $
and $\M_1(X^\nu)=\log 2$.  
\end{center}
 
On the other hand, the Aomoto complex for $\nu_\Z$ and $\nu_\K$ can be computed by the formula given in \cite[p.119]{Dim17}. A direct computation shows that 
\begin{center}
   $\tau_1(X,\nu_\Z)=4$ and
 $\beta_1(X,\nu_\K) =\begin{cases}
1, & \mathrm{if}\ \cha(\K)=2,\\
0 , & \mathrm{if}\ \cha(\K)\neq 2.
\end{cases} $ 
\end{center}
\eex

It is a long standing open question for hyperplane arrangement if $\sV^1_1(X,\C)$ is  combinatorially determined, see e.g. \cite[Problem 3.15]{Suc14}.  
 It is known that the  positive-dimensional components  of $\sV^1_1(X,\C)$ passing through the origin are combinatorially determined (see e.g. \cite[Theorem 6.1]{Dim17}) and they can be described by the multinet structure, for which we refer in particular to the work of Falk and Yuzvinsky \cite{FY}. 
  Since $\sV^1_1(X,\C)$ only depends on the the fundamental group of $X$, using Lefschetz hyperplane section theorem, one can assume that  $X$ is the complement of a line arrangement  in $\mathbb{CP}^2$.
Let us now recall the definition of multinet.
\bd \label{def multinet} Let $\sA$ be an arrangement in $\mathbb{CP}^2$.  A multinet on  $\sA$ is a partition of $\sA$ into $k\geq 3$ subsets $\sA_1,\cdots, \sA_k$, together with an assignment of multiplicities, $m\colon \sA \to \Z_{>0}$, and a subset $\mathfrak{X}$ of intersection points in $\sA$ satisfying the following conditions:
\begin{itemize}
\item[(a)] $\sum_{H\in \sA_i} m_H= \kappa$, independent of $i$;
\item[(b)] for each $H\in \sA_i$ and $H'\in \sA_j$ with $i\neq j$, $H\cap H' \in \mathfrak{X}$; 
\item[(c)] for each $x\in \mathfrak{X}$, the sum $\sum_{H\in \sA_i, x\in H} m_H$ is independent of $i$;
\item[(d)] for every $1\leq i \leq k$ and $H,H'\in \sA_i$, there is a sequence $\{H_0,H_1,\cdots,H_s\}$ with $H_0=H$ and $H_s=H'$ such that $H_{j-1}\cap H_j \notin \mathfrak{X}$ for all $1\leq j \leq s$.
\item[(e)] $\gcd(m_H\colon H\in \sA)$=1.
\end{itemize} Such multinet is called a $(k,\kappa)$-multinet.
\ed
Note that the multinet only depends on the intersection lattice of $\sA$.

By results of Pereira-Yuzvinsky \cite{PY} and Yuzvinsky \cite{Yuz}, $k=3$ or $4$ if $\vert \mathfrak{X}\vert>1$. Moreover, Yuzvinsky  conjectured that $k=4$ only happens for Hessian  arrangement of 12 planes in $\CP^2$, defined by the polynomial
\begin{center}
    $xyz \prod_{j=0}^2 \prod_{k=0}^2 (x+\omega^j y+\omega^k z )=0$ with $\omega=\exp{2\pi i/3}$.
\end{center}
As shown by Falk and Yuzvinsky \cite{FY}, and Pereira and Yuzvinsky \cite{PY}, every $(k,\kappa)$-multinet gives an  orbifold map $f\colon X \to C_{0,k}\subset \mathbb{CP}^1$ with $X=\mathbb{CP}^2-\sA$ by
$$f(x)\coloneqq [\prod_{H\in \sA_1} f_H^{m_H} (x), \prod_{H\in \sA_2} f_H^{m_H}(x)] .$$

\bd Let $X=\mathbb{CP}^2-\sA$ with $\sA$ being a essential line arrangement in $\mathbb{CP}^2$. We call a hyperbolic orbifold map $f\colon X\to C_{0,k}$ non-local, if it is not the restriction of another hyperbolic orbifold map over the complement of a sub-arrangement of $\sA$.
\ed 

\bt \cite{FY} Let $\sA$ be a essential hyperplane arrangement in $\mathbb{CP}^2$ with $X=\mathbb{CP}^2-\sA$.  Then every non-local hyperbolic orbifold map  $f\colon  X\to C_{0,k}$ with $k=3, 4$ is obtained from a $(k,\kappa)$-multinet for some $\kappa$.
\et

The above theorem shows that the existence of the non-local hyperbolic orbifold map  $f\colon  X\to C_{0,k}$ with $k=3$ or  $4$ are combinatorially determined (by the multinet structure). 
\begin{que} \label{que multi fiber} We have the following questions for $X=\CP^2\setminus \sA$ with $\sA$ being a essential line arrangement in $\CP^2$.
    \begin{itemize}
        \item[(i)] Is the existence of hyperbolic orbifold map  $f\colon  X\to C_{0,2}=\C^*$ combinatorially determined?
        \item[(ii)] Is the existence  of  non-local hyperbolic orbifold map  $f\colon  X\to C_{0,k}$ having  multiple fibers for $k=3$ or  $4$ combinatorially determined?
    \end{itemize}
\end{que}

 Assume that $\sA$ has $(d+1)$-many lines and let $f_i$ be the degree one defining equation for the $i$-th line $H_i\in \sA$ with  $0\leq i\leq d$. 
 \begin{assumption}\label{assumption tuple} For non-zero $a=(a_0,\cdots,a_d)\in \Z^{d+1}$, we assume that $\sum_{i=0}^d a_i=0$ and $\gcd(a_0,\cdots,a_d)=1$.
 \end{assumption} 
Then each $a\in \Z^{d+1}$ satisfying this assumption  gives a well-defined algebraic map $$f^a\coloneqq \prod_{i=0}^d f_i^{a_i}\colon X\to \C^*.$$
Note that $f^a$ has connected  generic fiber, since $\gcd(a_0,\cdots,a_d)=1$.  Then $f^a$ is a hyperbolic orbifold map if $f^a$ is surjective and has multiple fiber. 
In particular, every hyperbolic orbifold map has to be realized in this way. 
Note that $f^a$ induces a surjective map $$f^a_*\colon \pi_1(X)\to \pi_1(\C^*)\cong \Z. $$ This map only depends on $a$, which sends the meridian of $H_i$ to $a_i$. It gives a cohomology class in $H^1(X,\Z)$, say denoted by $\nu^a_\Z$. Proposition \ref{prop comb} shows that $f^a$ has multiple fibers only if $\tau_1(X,\nu^a_\Z)>1$,
which leads to the following combinatorics  question related to Question \ref{que multi fiber}(i).
\begin{que}
For which $a\in \Z^{d+1}$ satisfying Assumption \ref{assumption tuple}, we have  $\tau_1(X,\nu^a_\Z)>1$ ? 
\end{que}

On the other hand, we can also use Proposition \ref{prop comb} to study  Question \ref{que multi fiber}(ii).
Let $f\colon X\to C_{0,k}$ be a non-local orbifold map with $k=3$ or $4$. Then $f$ induces a surjective map  $$f_*\colon H_1(X,\Z)\to H_1(C_{0,k},\Z)\cong \Z^{k-1}. $$
Choosing any epimorphism $\nu\colon  \Z^{k-1} \to \Z,$ we obtain a composed epimorphism, still denoted by $\nu$. Denote corresponding cohomology class in $H^1(X,\Z)$  by $\nu_\Z$.  
Proposition \ref{prop comb} shows that $f$ has multiple fiber only if $\tau_1(X,\nu_\Z)>1$ for any such $\nu$.
Since $k=3$ or $4$, there are infinitely many choices for $\nu$, which gives a great possibility to show that $f$ has no multiple fiber in this case. We confirm this guess under some assumptions for multinet when $k=3$.

\begin{assumption} \label{assumption} Given a $(3,\kappa)$-multinet on  $\sA$ in $\CP^2$, we assume that
\begin{itemize}
\item[(a)] $m_H=1$ for any $H\in \sA_1$ or $ \sA_2$.
\item[(b)] If $H,H'\in \sA_1$ and $H\cap H' \notin \mathfrak{X}$, then there is no third hyperplane in $\sA$ passing through $H\cap H'$. So is $\sA_2$.
\item[(c)] Consider  $\sA_3$ and $\mathfrak{X}$. If there exists a point $x\in \mathfrak{X}$ such that $x$ is contained in precisely one hyperplane $H\in \sA_3$, one replaces $\sA_3$ and  $\mathfrak{X}$ by     deleting  the hyperplane $H$ and the point $x$, respectively.  We assume that one can run this procedure for $\sA_3$ and $\mathfrak{X}$ until there is no hyperplane left in $\sA_3$.
\end{itemize} 
\end{assumption}
Under this assumption, the author jointly with Li showed that the corresponding orbifold map $f$ has no multiple fibers. 
The technical conditions in Assumption \ref{assumption} are used to prove  $\tau_1(X,\nu_\Z)=1$ for some $\nu$.
\bp \cite[Proposition 5.9]{LL21} Let $f\colon X\to C_{0,3}$ be an orbifold map corresponding to  a $(3,\kappa)$-multinet structure for $\sA$ in $\CP^2$ with $X=\CP^2 \setminus \sA$. If this multinet satisfies Assumption \ref{assumption}, $f$ has no multiple fibers.
\ep 

\bex We give a list of examples which satisfies Assumption \ref{assumption}:
 Figure 1(a), Figure 2, Figure 3(a),(b) in \cite{FY} with $A,B,C$ corresponding to $\sA_1,\sA_2,\sA_3$ and Figure 1 in  \cite{Suc14}.
It is also easy to find examples which do not satisfy Assumption \ref{assumption}, see Figure 1(b) in \cite{FY} or the monomial (alias Ceva)  arrangement (with $m\geq 3$) in \cite[Example 6.11]{Dim17}.  
\eex

Consider a non-local  orbifold map $f\colon X \to C_{0,k} $ with $k=2,3,4$. 
When $k=2$, there are plenty of examples where $f$ has multiple fibers, e.g. see Example \ref{ex 2}. 
In fact, Suciu discovered a combinatorial construction to produce such orbifold maps having multiple fibers from a given multinet,  
called the pointed multinet structure, see e.g. \cite[Section 3.6]{Suc14}. He asked if all such orbifold maps appear in this way, see \cite[Problem 3.12]{Suc14}.

On the other hand, to the author's knowledge, 
 no examples with $f$ having multiple fibers are known when $k=3,4$. 
\begin{que}
    Is it possible to show that $f$ always has no multiple fibers when $k=3$ or $k=4$?
\end{que}
All the questions in this section should be considered as sub-questions of the aforementioned question: 
\begin{que}
    Is $\sV^1_1(X,\C)$ combinatorially determined for hyperplane arrangement complement $X$?
\end{que} 

\section{Proof of results}

\subsection{Proof of Proposition \ref{prop lim}}
Set  $\Z[\T^\pm]=\Z[\Z^n].$
Let $M$ be a finitely generated $\Z[\T^\pm]$-module. Consider  a presentation of $M$ by an exact sequence 
$$ \Z[\T^\pm]^{b'} \overset{\partial}{\to} \Z[\T^\pm]^b \to M \to 0$$
where $\partial$ is given by a matrix of size $b'\times b$ with entries in $\Z[\T^\pm]$.
For any $j\geq 0$, let $\Delta^j(M)$ denote the the greatest
common divisor of all the $(b -j)\times (b-j)$ minors of $\partial$ with the convention that $\Delta^j(M)=1$  if $j\geq b$. Clearly $\Delta^j(M)$ divides $\Delta^{j-1}(M)$ for every $j\geq 1$.  It 
is known that $\Delta^j(M)$ only depends on $M$, not on any particular presentation. The rank of the module $M$ over $\Z[\T^\pm]$ is the dimension of the vector space $M\otimes_{\Z[\T^\pm]} \mathcal{F}'$ 
over the fractional field $\mathcal{F}'$ of $\Z[\T^\pm]$.  If M has rank $r$, then $\Delta^j (M) = 0$ if $j <r$, and $\Delta^{j}(M)=\Delta^{j-r}(\Tor( M))$, where $$\Tor(M)=\{ x\in M \mid fx=0 \text{ for some non-zero } f\in \Z[\T^\pm]\}.$$

As in  \cite[section 1.3]{Le}, a $\Z[\T^\pm]$-module $M$ is called pseudo-zero if for every prime ideal $P$ of height 1, the localization $M_P$ is 0.  A $\Z[\T^\pm]$-module morphism $M_1 \to M_2$ is a pseudo-isomorphism if the kernel and co-kernel
are pseudo-zero. Moreover, two finitely generated torsion $\Z[\T^\pm]$-modules $M_1$ and $M_2$ are pseudo-isomorphic if and
only if $\Delta^j(M_1)=\Delta^j(M_2)$ for any  $j\geq 0$.

To prove Proposition \ref{prop lim}, let us first recall a theorem proved by Le. 

\bt \cite[Theorem 6]{Le} \label{thm pseudo} Assume that $M_1$ and $M_2$ are two pseudo-isomorphic finitely generated $\Z[\T^\pm]$ modules. Then $\vert \tor_\Z \big( M_1\otimes_{\Z[\T^\pm]} \Z[\Z^n/\Gamma]\big)\vert $ and   $\vert \tor_\Z \big(M_2\otimes_{\Z[\T^\pm]} \Z[\Z^n/\Gamma]\big)\vert $ have the same growth rate in the sense that
$$ \lim_{\langle\Gamma\rangle\rightarrow\infty}\big(\frac{\log  \vert \tor_\Z \big(M_1\otimes_{\Z[\T^\pm]} \Z[\Z^n/\Gamma]\big)\vert }{\vert  \Z^n/\Gamma\vert}- \frac{\log  \vert \tor_\Z \big(M_2\otimes_{\Z[\T^\pm]} \Z[\Z^n/\Gamma]\big)\vert }{\vert  \Z^n/\Gamma\vert}\big)=0. $$
Here $\vert \tor_\Z \big(M_i\otimes_{\Z[\T^\pm]} \Z[\Z^n/\Gamma]\big) \vert$ denotes the order of the torsion part of $ M_i\otimes_{\Z[t^\pm]} \Z[\Z^n/\Gamma]$ as an abelian group.
\et 
Note that the limit in this theorem  is the ordinary limit, not the upper limit. The following result is  a direct application of this theorem and Theorem \ref{thm cyclotomic}. 

\bp \label{prop limit} Let $M$ be a finitely generated $\Z[\T^\pm]$-modules with rank $r$. If $\Delta^r(M)$ is  a product of some generalized cyclotomic polynomials with a non-zero integer $c$,
then $$\lim_{\langle\Gamma\rangle\rightarrow\infty} \frac{\log  \vert \tor_\Z \big(M\otimes_{\Z[\T^\pm]} \Z[\Z^n/\Gamma]\big)\vert }{\vert  \Z^n/\Gamma\vert}=\log \vert c \vert.$$ 
\ep 
\begin{proof}
Consider the sub-module $\Tor(M)$ of $M$. We have that $$\Delta^j(M)=\Delta^{j-r}(\Tor(M)).$$ Moreover, it follows from Theorem \ref{thm cyclotomic} and  \cite[Theorem 1.2, Proposirion 3.4]{Le} that
$$\lim_{\langle\Gamma\rangle\rightarrow\infty}\big(\frac{\log  \vert \tor_\Z \big(M\otimes_{\Z[\T^\pm]} \Z[\Z^n/\Gamma]\big)\vert }{\vert  \Z^n/\Gamma\vert}- \frac{\log  \vert \tor_\Z \big(\Tor(M)\otimes_{\Z[\T^\pm]} \Z[\Z^n/\Gamma]\big)\vert }{\vert  \Z^n/\Gamma\vert}\big)=0.$$

So without loss of generality, we assume that $M$ is a finitely generated torsion $\Z[\T^\pm]$-module.
Assume that $\Delta^j(M)=1$ for $j>k$. 
    Consider a finitely generated torsion $\Z[\T^\pm]$-module  $$M'\coloneqq \bigoplus_{j=0}^k \frac{ \Z[\T^\pm]}{\big(\Delta^j(M)/\Delta^{j+1}(M) \big)} .$$ It is easy to see that $\Delta^j(M)=\Delta^j(M')$ for any $j\geq 0$. Then $M$ is pseudo-isomorphic to $M'$ by \cite[Theorem 3.5]{Hil}. Using Theorem \ref{thm pseudo}, the proof is reduced to $M'$. Note that $\Delta^{j+1}(M)$ divides $\Delta^{j}(M)$ for all $j\geq 0$, hence $ \Delta^j(M)/\Delta^{j+1}(M)$ is a also product of some generalized cyclotomic polynomials and some non-zero integer.    Since $M'$ is a direct sum, the proof is further reduced to the module with the type $N=\Z[\T^\pm]/(h)$, where $h$ a product of some generalized cyclotomic polynomials and some non-zero integer $d$.

Set $h'=\frac{h}{d}$ and $N'=\Z[\T^\pm]/(h')$. 
    Clearly we have $$\vert \tor_\Z\big( N\otimes_{\Z[\T^\pm]} \Z[\Z^n/\Gamma]\big)\vert= \vert d \vert^{\vert \Z^n/\Gamma\vert} \cdot \vert \tor_\Z \big(N'\otimes_{\Z[\T^\pm]} \Z[\Z^n/\Gamma]\big)\vert.$$
    Since $h'$ is  a product of some generalized cyclotomic polynomials, 
    using Theorem \ref{thm cyclotomic} and Theorem \cite[Theorem 3]{Le} we have 
    $$\lim_{\langle\Gamma\rangle\rightarrow\infty} \frac{\log  \vert \tor_\Z \big( N'\otimes_{\Z[\T^\pm]} \Z[\Z^n/\Gamma]\big)\vert }{\vert  \Z^n/\Gamma\vert}=0.$$
 Then the claim follows.   
\end{proof}
\begin{proof}[Proof of Proposition \ref{prop lim}]

As shown in the proof of \cite[Theorem 4]{Le}, we have $$ |H_{i}(X^{\nu,\Gamma},\mathbb{Z})_{\mathrm{tor}}|=\vert \tor_\Z \big(M\otimes_{\Z[t^\pm]} \Z[\Z^n/\Gamma]\big)\vert $$ with $M=\coker \partial_i$, see \cite[(54) on page 745]{Le}. 
Clearly $M$ has a presentation  by the exact sequence 
$$ C_{i+1}(X^\nu,\Z) \overset{\partial_i}{\to} C_{i}(X^\nu,\Z) \to M\to 0$$
and $\rank M=\rank  C_{i}(X^\nu,\Z)-\rank \partial_i$. Then we have $\Delta^{\rank M}(M)=\Delta_i(X^\nu)$ by definition.  Therefore the claim follows from Proposition \ref{prop limit}.
\end{proof}

\subsection{Proof of Proposition \ref{prop betti number}}
  The first equality follows from the proof of \cite[Propsition 2.2]{LL23}. The second equality follows directly from the fact that  computing $\rank H_i(X^\nu,\K)$ over the fraction field is same as computing the corresponding homology over a general point.

  The second claim follows from  the fact that every irreducible component of the zero locus of $\Delta_i(X^\nu)$ in $(\C^*)^n$ gives a prime ideal with height one in $\C[\T^\pm]$. Taking localization of $\C[\T^\pm]$ over this primed ideal, one gets a PID, say denoted by $R$. 
  Consider the localized complex $C_*(X^\nu,\C)\otimes_{\C[\T^\pm]} R$. 
  Then the second claim follows easily by analysing this localized complex of modules over PID, see e.g. \cite[Section 2.2]{BLW}.

\subsection{Proof of Proposition \ref{prop lim hodge}}
Let $\Pic^0(X)$  denote the identity component of the Picard variety of $X$, which parameterise topologically trivial line bundles on $X$. Set $$S^{p,q}_k(X)\coloneqq \{\sL\in\Pic^0(X)\mid \dim H^q(X, \Omega^p_X \otimes \sL) \geq k \}.$$
Green and Lazarsfeld showed \cite{GL91} that $S^{p,q}_k(X)$ is a finite union of torsion translated sub-abelian variety of $\Pic^0(X)$. In particular, these are sub-varieties of $\Pic^0(X)$.

Consider the moduli space of unitary rank one $\C$-local system on $X$ defined as
$$ \mathrm{Uchar}^0(X)\coloneqq \homo (H_1(X,\Z)/\{\mathrm{torsion} \}, U(1)) .$$
Any unitary local system $L\in \mathrm{Uchar}^0(X)$ corresponds to a topologically trivial line bundle $\sL=L\otimes_\C \mathcal{O}_X$. This gives a one-to-one correspondence between $\mathrm{Uchar}^0(X)$ and $\Pic^0(X)$ as real torus.
By slight modification of the usual Hodge theoretic arguments, one gets the following Hodge decomposition
\be \label{Hodge decomposition}   H^i(X,L) \cong \bigoplus_{p+q=i} H^q(X,\Omega^p_X \otimes \sL).
\ee

The epimorphism $\nu$ induces an embedding  $\nu^\#\colon (U(1))^n \hookrightarrow \mathrm{Uchar}^0(X).$ Let $\im \nu^\#$ denote the corresponding image in $\Pic^0(X)$ by the above identification.
Since $S^{p,q}_k(X)$ are sub-varieties of $\Pic^0(X)$, we can pick   a generic unitary local system $L'\in \im \nu^\#$ with $\sL'=L'\otimes \sO_X$ (to be precise, $L'$ is generic with respect to all degrees $(p,q)$ with $0\leq p,q \leq \dim X$). By semi-continuity, for any unitary local system $L$ and its corresponding line bundle $\sL$ we have 
$$ \dim H^q(X, \Omega^p_X \otimes \sL') \leq \dim H^q(X, \Omega^p_X \otimes \sL).$$

 Let $\pi\colon X^{\nu,\Gamma} \to X$ denote the covering map. 
Then $$\pi_* \sO_X=\bigoplus \sL_\lambda$$
where $\sL_\lambda=L_\lambda\otimes \sO_X$ and the direct sum runs over all unitary local system $L_\lambda$ with $\lambda\in \homo(\Z^n/\Gamma,\C^*)$.
 We have 
$$ \dim H^q(X^\Gamma, \Omega^p_{X^\Gamma}) =\sum_\lambda \dim H^q(X, \Omega^p_X \otimes \sL_\lambda)\geq \vert\Z^n/\Gamma\vert \cdot \dim H^q(X, \Omega^p_X \otimes \sL').$$
Meanwhile, the Hodge decomposition  gives us 
$$\sum_{p+q=i}\dim H^q(X^\Gamma, \Omega^p_{X^\Gamma})=\dim H^i(X^\Gamma,\C).$$
Putting all together, we have that 
$$\sum_{p+q=i}\dim H^q(X, \Omega^p_X \otimes \sL') \leq  \sum_{p+q=i}\dfrac{ \dim H^q(X^\Gamma, \Omega^p_{X^\Gamma})}{\vert\Z^n/\Gamma\vert} =\dfrac{ \dim H^i(X^\Gamma,\C)}{\vert\Z^n/\Gamma\vert} .$$

$\sV^i_k(X)$ are also closed sub-varieties.  By further requiring the unitary local system $L'$ generic with respect to the jump loci, Proposition \ref{prop betti number} and Universal Coefficients Theorem imply that the right hand of the above inequality converges to $\dim H^i(X,L')$ when $\langle \Gamma \rangle \to \infty$.
Note that  the Hodge decomposition (\ref{Hodge decomposition}) gives us $$ \sum_{p+q=i}\dim H^q(X, \Omega^p_X \otimes \sL') =\dim H^i(X,L').$$
Therefore we have $$\lim_{\langle\Gamma\rangle\rightarrow\infty}\frac{h^{p,q}(X^{\nu,\Gamma})}{|\Z^n/\Gamma|} =\dim H^q(X, \Omega^p_X \otimes \sL'). $$

\bibliographystyle{amsalpha}

\end{document}